\documentclass[10pt]{article}

\usepackage{amsthm}
\usepackage{amsmath}
\usepackage{amsfonts}
\usepackage{amssymb}
\usepackage{mathrsfs}
\usepackage{stmaryrd}
\usepackage{dsfont}
\usepackage{enumerate}
\usepackage{soul}
\usepackage[colorlinks=true,urlcolor=blue]{hyperref}
\usepackage[applemac]{inputenc}

\newcommand{\N}{\mathbb{N}}

\newcommand{\R}{\mathbb{R}}

\newcommand{\dom}{\mbox{\rm dom}\,}
\renewcommand{\O}{\Omega}
\renewcommand{\o}{\omega}

\newcommand{\vE}{{\cal E}}
\newcommand{\vF}{{\cal F}}

\newcommand{\vO}{{\cal O}}

\newcommand{\Mr}{\mathscr{M}}
\newcommand{\vphi}{\varphi}
\newcommand{\eps}{\varepsilon}
\newcommand{\la}{\lambda}

\newcommand{\ovl}{\overline}
\newcommand{\udl}{\underline}
\newcommand{\vlim}{\lim\limits}

\newcommand{\vmax}{\max\limits}

\newcommand{\vsup}{\sup\limits}

\newcommand{\vint}{\int\limits}

\newcommand{\tends}{\longrightarrow}

\newcommand{\wt}{\widetilde}
\newcommand{\wh}{\widehat}

\newcommand{\loc}{\mathrm{loc}}

\renewcommand{\b}{\mathrm{b}}
\renewcommand{\d}{\mathrm{d}}

\newcommand{\dist}{\mathrm{dist}}

\renewcommand{\le}{\leqslant}
\renewcommand{\ge}{\geqslant}

\newcommand{\bs}{\boldsymbol}

\newcommand{\p}{\prime}
\newcommand{\pp}{{\prime\prime}}

\DeclareMathOperator{\supp}{supp}

\DeclareMathOperator{\sign}{sign}

\DeclareMathOperator{\crit}{crit}
\DeclareMathOperator{\argmin}{argmin}
\DeclareMathOperator{\inte}{int}

\numberwithin{equation}{section}

\newtheorem{thm}{Theorem}[section]
\newtheorem{prop}[thm]{Proposition}
\newtheorem{cor}[thm]{Corollary}
\newtheorem{lem}[thm]{Lemma}

\theoremstyle{definition}
\newtheorem{rmk}[thm]{Remark}
\newtheorem{defi}[thm]{Definition}
\newtheorem{exa}[thm]{Example}

\newenvironment{proof*}{\noindent{\bf Proof.}}{\qed}
\newenvironment{vproof}[1]{\noindent{\bf Proof #1}}{\qed}

\title{On damped second-order gradient systems}
\author{\sc Pascal Bégout\footnote{TSE (Institut de Mathématiques de Toulouse, Université Toulouse I Capitole), Manufacture des Tabacs, 21 allée de Brienne, 31015 Toulouse, Cedex 06, France}, Jérôme Bolte\footnote{TSE (GREMAQ, Université Toulouse I Capitole), Manufacture des Tabacs, 21 allée de Brienne, 31015 Toulouse, Cedex 06, France. Effort sponsored by the Air Force Office of Scientific Research, Air Force Material Command, USAF, under grant number FA9550-14-1-0056. This research also benefited from the support of the ``FMJH Program Gaspard Monge in optimization and operations research''} \,and Mohamed Ali Jendoubi\footnote{Université de Carthage, Institut préparatoire aux études scientifiques et techniques, BP 51, 2080 La Marsa, Tunisia}}
\date{}

\begin{document}

\maketitle

\vspace{-0.5cm}

\begin{abstract}
Using small deformations of the total energy, as introduced in \cite{MR1616968}, we establish that damped second order gradient systems 
\begin{gather*}
u^\pp(t)+\gamma u^\p(t)+\nabla G(u(t))=0,
\end{gather*}
may be viewed as quasi-gradient systems. In order to study the asymptotic behavior of these systems, we prove that any (nontrivial) desingularizing function appearing in KL inequality satisfies $\vphi(s)\ge c\sqrt s$ whenever the original function is definable and $C^2.$ Variants to this result are given. These facts are used in turn to prove that a desingularizing function of the potential $G$ also desingularizes the total energy and its deformed versions. Our approach brings forward several results interesting for their own sake: we provide an asymptotic alternative for quasi-gradient systems, either a trajectory converges, or its norm tends to infinity. The convergence rates are also analyzed by an original method based on a one-dimensional worst-case gradient system.

We conclude by establishing the convergence of solutions of damped second order systems in various cases including the definable case. The real-analytic case is recovered and some results concerning convex functions are also derived.
\end{abstract}

{\let\thefootnote\relax\footnotetext{E-mail:
\htmladdnormallink{{\footnotesize\udl{\texttt{Pascal.Begout@tse-fr.eu}}}}{mailto:Pascal.Begout@math.cnrs.fr}($^\ast$),
\htmladdnormallink{{\footnotesize\udl{\texttt{Jerome.Bolte@tse-fr.eu}}}}{mailto:Jerome.Bolte@tse-fr.eu}$(^\dagger),$
\htmladdnormallink{{\footnotesize\udl{\texttt{ma.jendoubi@fsb.rnu.tn}}}}{mailto:ma.jendoubi@fsb.rnu.tn}$(^\ddagger)$
}}
{\let\thefootnote\relax}
{\let\thefootnote\relax}

\vspace{-0.5cm}

\tableofcontents

\section{Introduction}
\label{intro}

\subsection{A global view on previous results}
In this paper, we develop some new tools for the asymptotic behavior as $t$ goes to infinity of solutions $u:\R_+\tends\R^N$ of the following second order system
\begin{gather}
\label{second1}
u^\pp(t) + \gamma u^\p (t) + \nabla G(u(t)) = 0, \quad t\in\R_+.
\end{gather}
Here, $\gamma>0$ is a positive real number which can be seen as a {\em damping coefficient}, $N\ge1$ is an integer and $G\in C^2(\R^N)$ is a real-valued function. In Mechanics, \eqref{second1} models, among other problems, the motion of an object subject to a force deriving from a potential $G$ (e.g. gravity) and to a viscous friction force $-\gamma u^\p.$ In particular, the above may be seen as a qualitative model for the motion of a material point subject to gravity, constrained to evolve on the graph of $G$ and subject to a damping force, further insights and results on this view may be found in \cite{MR1909449,MR1918942}. This type of dynamical system has been the subject of several works in various fields and along different perspectives, one can quote for instance \cite{MR1753136} for Nonsmooth Mechanics, \cite{MR3218822,MR3500980} for recent advances in Optimization and \cite{zbMATH03238176} for pioneer works on the topic, partial differential equations and related aspects \cite{MR2328934,MR2018329,MR3193987}.
\medskip \\
The aim of this work is to provide a deeper understanding of the asymptotic behavior of such a system and of the mechanisms behind the stabilization of trajectories at infinity (making each bounded orbit approach some specific critical point). Such behaviors have been widely investigated for gradient systems,
\begin{gather*}
u^\p(t)+\nabla G(u(t))=0,
\end{gather*}
for a long time now. The first decisive steps were made by {\L}ojasiewicz for analytic functions through the introduction of the so-called gradient inequality \cite{MR0160856,Lo65}. Many other works followed among which two important contributions: \cite{MR0377609} for convex functions and \cite{MR1644089} for definable functions. Surprisingly the asymptotic behavior of the companion dynamics \eqref{second1} has only been ``recently" analyzed. The motivation for studying \eqref{second1} seems to come from three distinct fields PDEs, Mechanics and Optimization. Out of the convex realm \cite{MR1104346,MR1760062}, the seminal paper is probably \cite{MR1616968}. Like many of the works on gradient systems the main assumption, borrowed from {\L}ojasiewicz original contributions, is the analyticity of the function -- or more precisely the fact that the function satisfies the {\L}ojasiewicz inequality. This work paved the way for many developments: convergence rates studies \cite{MR1829143}, extension to partial differential equations \cite{MR727703,MR1609269,MR1616964,MR1714129,MR1680877,MR2328934,MR1986700,MR1800136,MR2041509,MR2772353,MR3193987}, use of various kind of dampings \cite{MR2429439,MR2511558} (see also \cite{MR2852207,MR3093225,MR3189085,MR2727316}). Despite the huge amount of subsequent works, some deep questions remained somehow unanswered; in particular it is not clear to see:
\begin{itemize}
\item[--]
{\em What are the exact connections between gradient systems and damped second-order gradient systems?}
\item[--]
{\em Within these relationships, how central is the role of the properties/geometry of the potential function $G$?}
\end{itemize}
Before trying to provide some answers, we recall some fundamental notions related to these questions; they will also constitute the main ingredients in our analysis of~\eqref{second1}.

\noindent
\textbf{Quasi-gradient fields.} The notion is natural and simple: a vector field $V$ is called \textit{quasi-gradient} for a function $L$ if it has the same singular point (as $\nabla L$) and if the angle $\alpha$ between the field $V$ and the gradient $\nabla L$ remains acute and bounded away from $\pi/2.$ Proper definitions are recalled in Section~\ref{Lyapunov}. Of course, such systems have a behavior which is very similar to those of gradient systems (see Theorem~\ref{thmsqg}). We refer to \cite{MR2901252} and the references therein for further geometrical insights on the topic.

\noindent
\textbf{Liapunov functions for damped second order gradient systems.} The most striking common point between \eqref{second1} and gradient systems is that of a ``natural" Liapunov function. In our case, it is given by the {\em total energy}, sum of the potential energy and the kinetic energy,
\begin{gather*}
E_T(u,v)=G(u)+\frac{1}{2}\|v\|^2.
\end{gather*}
The above is a Liapunov function in the phase space, more concretely
\begin{eqnarray*}
\frac{\d}{\d t}E_T(u(t),u^\p(t))	&	=	& \frac{\d}{\d t}\Big(\frac12|u^\p(t)|^2+G(u(t))\Big) \\
						&	=	& -\gamma \|u^\p(t)\|^2.
\end{eqnarray*} Contrary to what happens for classical gradient systems the vector field associated with \eqref{second1} is not strictly Lyapunov for $E_T$: it obviously degenerates on the subspace $[v=0]$ (or $[u^\p=0]$). The use of $E_T$ is however at the heart of most results attached to this dynamical system.

\noindent
\textbf{KL functions.} A KL function is a function whose values can be reparametrized in the neighborhood of each of its critical point so that the resulting functions become {\em sharp}$^($\footnote{That is, the norms of its gradient remain bounded away from zero.}$^).$ More formally, $G$ is called KL on the slice of level lines $[0<G<r_0]\stackrel{\text{def}}{=}\big\{u\in \R^N;0<G(u)<r_0\big\},$ if there exists $\vphi\in C^0\big([0,r_0)\big)\cap C^1(0,r_0)$ concave such that $\vphi(0)=0,$ $\vphi^\p>0$ and
\begin{gather*}
\|\nabla (\vphi\circ G)(u)\|\ge 1, \quad \forall u\in[0<G<r_0].
\end{gather*}

\noindent
Proper definitions and local versions can be found in the next section. The above definition originates in \cite{MR2592958} and is based on the fundamental work of Kurdyka \cite{MR1644089}, where it was introduced in the framework of o-minimal structure$^($\footnote{A far reaching concept that generalizes semi-algebraic or (globally) subanalytic classes of sets and functions.}$^)$ as a generalization of the famous {\L}ojasiewicz inequality. 

\noindent
KL functions are central in the analysis of gradient systems, the readers are referred to \cite{MR2592958} and the references therein.

\noindent
\textbf{Desingularizing functions.} The function appearing above, namely $\vphi,$ is called a {\em desingularizing function}: the faster $\vphi^\p$ tends to infinity at $0,$ the flatter is $G$ around critical points. As opposed to the {\L}ojasiewicz gradient inequality, this behavior, in the o-minimal world, is not necessarily of a ``power-type". Highly degenerate functions can be met, like for instance $G(u)=\exp\big(-1/p^2(u)\big)$ where $p:\R^N\tends\R$ is any real polynomial function. This class of functions belongs to the log-exp structure, an o-minimal class that contains semi-algebraic sets and the graph of the exponential function \cite{MR1398816}. Finally, observe that if it is obvious that $\vphi$ might have an arbitrarily brutal behavior at $0,$ it is also pretty clear that the smoothness of $G$ is related to a lower-control of the behavior of $\vphi,$ for instance we must have $\vphi^\p(0)=\infty$ -- which is not the case in general in the nonsmooth world (see e.g. \cite{MR2338451}). 

\subsection{Main results}

Several auxiliary theorems were necessary to establish our main result, we believe they are interesting for their own sake. Here they are:

-- An asymptotic alternative for quasi-gradient systems: either a trajectory converges or it escapes to infinity,

-- A general convergence rate result for the solutions of the gradient systems that brings forward a worst-case gradient dynamical system in dimension one,

-- Lower bounds for desingularizing functions of $C^2$ KL functions.

\bigskip

\noindent
We are now in position to describe the strategy we followed in that paper for the asymptotic study of the damped second order gradient system \eqref{second1}. Our method was naturally inspired by the Liapunov function provided in \cite{MR1616968}. 

\begin{enumerate}
\item
First we show that $E_T$ can be slightly and ``semi-algebraically" (respectively, definably) deformed into a smooth function $E^\mathrm{def}_T,$ so that the gradient of the new energy $\nabla E^\mathrm{def}_T$ makes an uniformly acute angle with the vector field associated with \eqref{second1} -- this property only holds on bounded sets of the phase space.
The system \eqref{second1} appears therefore as a quasi-gradient system for $E^\mathrm{def}_T.$

\item
In a second step we establish/verify that the solutions of the quasi-gradient systems converge whenever they originate from a KL function. 

We also provide rates of convergence and we explain how they may be naturally and systematically derived from a one-dimensional worst-case gradient dynamics.

At this stage it is possible to proceed abstractly to the proof of the convergence of solutions to \eqref{second1} in several cases. For instance the definable case: we simply have to use the fact that $E^\mathrm{def}_T$ is definable whenever $G$ is, so it is a KL function and the conclusion follows.

Although direct and fast, this approach has an important drawback from a conceptual viewpoint since it relies on a desingularizing function attached to an auxiliary function $E^\mathrm{def}_T$ whose meaning is unclear. Whatever perspectives we may adopt (Mechanics, Optimization, PDEs), an important question is indeed to {\em understand what happens when $G$ is KL and how the desingularizing function of $G$ actually impacts the convergence of solutions to \eqref{second1}.}

\item
We answer to this question in the following way.
\begin{enumerate}
\item
We prove that desingularizing functions of $C^2$ definable functions have a lower bound. Roughly speaking, we prove that for nontrivial critical points the desingularizing function has the property $\vphi (s)\ge c\sqrt{s}$ \Big(or
equivalently$^($\footnote{Recall that $\vphi$ is definable.}$^)$ $\vphi^\p(s)\ge\left.\frac{c^\p}{\sqrt{s}}\right).$
\item
We establish that if $\vphi$ is definable and desingularizing for $G$ at $\ovl u$ then it is desingularizing for both $E_T$ and $E^\mathrm{def}_T$ at $(\ovl u, 0).$
\end{enumerate}
\item
We conclude by combining previous results to obtain in particular the convergence of solutions to \eqref{second1} under definability assumptions. We also provide convergence rates that depend on the desingularizing function of $G,$ {\em i.e. on the geometry of the potential.}
\end{enumerate}

\noindent
We would like to point out and emphasize two facts that we think are of interest. First the property $\vphi (s)\ge c\sqrt{s}$ (see Lemma~\ref{lower} below) is a new result and despite its ``intuitive" aspect the proof is nontrivial. We believe it has an interest in its own sake.

\medskip
\noindent
More related to our work is the fact that (in the definable case and in many other relevant cases) our results show that the desingularizing function of $G$ is conditioning the asymptotic behavior of solutions of the system. Within an Optimization perspective this means that the ``complexity", or at least the convergence rate, of the dynamical system is entirely embodied in $G$ when $G$ is smooth. From a mechanical viewpoint, stabilization at infinity is determined by the conditioning of $G$ provided the latter is smooth enough; in other words the intuition that for large time behaviors, the potential has a predominant effect on the system is correct -- a fact which is of course related to the dissipation of the kinetic energy at a ``constant rate".

\medskip
\noindent
\textbf{Notation.} The finite-dimensional space $\R^N$ $(N\ge1)$ is endowed with the canonical scalar product $\langle\:.\:,\:.\:\rangle$ whose norm is denoted by $\|\:.\:\|.$ The product space $\R^N\times\R^N$ is endowed with the natural product metric which we still denote by $\langle\:.\:,\:.\:\rangle.$ We also define for any $\ovl u\in\R^N$ and $r>0,$ $B(\ovl u,r)=\{u\in\R^N; \|u-\ovl u\|<r\}.$ When $S$ is a subset of $\R^N$ its interior is denoted by $\inte S$ and its closure by $\ovl S.$ If $F:\R^N\tends\R$ is a differentiable function, its gradient is denoted by $\nabla F.$ When $F$ is a twice differentiable function, its Hessian is denoted by $\nabla^2F.$ The set of critical points of $F$ is defined by
\begin{gather*}
\crit F=\Big\{u\in \R^N;\nabla F(u)=0\Big\}.
\end{gather*}
This paper is organized as follows. In Section~\ref{S:KL_ineq}, we provide a lower bound for desingularizing function of $C^2$ functions under various assumptions, like definability (Proposition~\ref{propG} and Lemma~\ref{lower}). In Section~\ref{quasi_gradient}, we recall the behavior of a first order system having a quasi-gradient structure for some KL function and we provide an asymptotic alternative (Theorem~\ref{thmsqg}). In Theorem~\ref{worstquas}, the convergence rate of any solution to a first order system having a quasi-gradient structure is proved to be better than that of a one-dimensional worst-case gradient dynamics (various known results are recovered in a transparent way). Finally, we establish that any function which desingularizes $G$ in \eqref{second1} also desingularizes the total energy and various relevant 
deformation of the latter (Proposition~\ref{G_Energy}). In Section~\ref{convergence}, we study the asymptotic behavior of solutions to \eqref{second1} (Theorem~\ref{thmmain}) while in Section~\ref{consequences}, we describe several consequences of our main results. Appendix provides, for the comfort of the reader, some elementary facts on o-minimal structures.

\section{Structural results: lower bounds for desingularizing functions of $\bs{C^2}$ functions}
\label{S:KL_ineq}

To keep the reading smooth and easy, we will not formally define here o-minimal structure. The definition is postponed in Appendix. Let us however recall, at this stage, that the simplest o-minimal structure (containing the graph of the real product) is given by the class of real semi-algebraic sets and functions. A semi-algebraic set is the finite union of sets of the form
\begin{gather}
\label{SA}
\Big\{u\in\R^N; \; p(u)=0, \: p_i(u)<0, \forall i\in I\Big\},
\end{gather}
where $I$ is a finite set and $p,\{p_i\}_{i\in I}$ are real polynomial functions.
\bigskip \\
Let us recall a fundamental concept for dissipative dynamical systems of gradient type.

\begin{defi}[\textbf{Kurdyka-{\L}ojasiewicz property and desingularizing function}]
\label{KL_property}$\;$\\
Let $G:\R^N\tends\R$ be a differentiable function.
\begin{itemize}
\item[(i)]
We shall say that $G$ has the \textit{KL property} at $\ovl u\in\R^N$ if there exist $r_0>0,$ $\eta>0$ and
$\vphi\in C([0,r_0);\R_+)$ such that
\begin{enumerate}
\item
$\vphi (0)=0,$ $\vphi\in C^1((0,r_0);\R_+)$ concave and $\vphi^\p$ positive on $(0,r_0),$
\item
$u\in B(\ovl u,\eta)\implies | G(u)-G(\ovl u)|<r_0$; and for each $u\in B\big(\ovl u,\eta\big),$ such that $G(u)\neq G(\ovl u),$\begin{gather}
\label{Kurd}
\big\|\nabla(\vphi\circ |G(\:.\:)-G(\ovl u)|)(u) \big\|\ge1.
\end{gather}
\end{enumerate}
Such a function $\vphi$ is called a \textit{desingularizing function} of $G$ at $\ovl u$ on $B(\ovl u,\eta).$
\item[(ii)] The function $G$ is called a {\em KL function} if it has the KL property at each of its points.
\end{itemize}
\end{defi}

\noindent
The following result is due to {\L}ojasiewicz in its real-analytic version (see e.g.~\cite{Lo65,MR0160856}), it was generalized to o-minimal structures and considerably simplified by Kurdyka in \cite{MR1644089} (see Appendix). 

\begin{thm}[\textbf{Kurdyka-{\L}ojasiewicz inequality \cite{MR1644089}}$^{\bs(}$\footnote{See comments in Appendix.}$^{\bs)}$]
\label{KL}
Let $\vO$ be an o-minimal structure and let $G\in C^1(\R^N;\R)$ be a definable function. Then $G$ is a KL function.
\end{thm}

\begin{rmk}
\label{rmkKL} 
(a) Theorem~\ref{KL} is of course trivial when $\ovl u\not\in\crit G$ -- take indeed, $\vphi(s)=cs$ where $c=\frac{1+\eps}{\|\nabla G(\ovl u)\|}$ and~$\eps>0.$ \\
(b) Restrictions of real-analytic functions to compact sets included in their (open) domain belong to the o-minimal structure of globally analytic sets \cite{MR1404337}. They are therefore KL functions (see indeed Example~\ref{exomin}). In some o-minimal structures there are nontrivial functions for which all derivatives vanish on some nonempty set, like $G(u)=\exp(-1/f^2(u))$ where $f\neq0$ is any smooth semi-algebraic function achieving the value $0^($\footnote{This function is definable in the log-exp structure of Wilkie \cite{MR1398816}.}$^)$ (see also Example~\ref{exomin}). For these cases, $\varphi$ is not of power-type -- as it is the case when $G$ is semi-algebraic or real-analytic. Other types of functions satisfying the KL property in various contexts are provided in~\cite{MR2674728} (see also Corollary~\ref{convex}). \\
(c) Desingularizing functions of definable functions can be chosen to be definable, strictly concave and $C^k$ (where $k$ is arbitrary).
\end{rmk}

\noindent
The following trivial notion is quite convenient. 

\begin{defi}[\textbf{Trivial critical points}]
\label{nonflatcrit}
A critical point $u$ of a differentiable function $G:\R^N\tends\R$ is called \textit{trivial} if $u\in\inte\crit G.$ It is \textit{nontrivial} otherwise. Observe that $u$ is nontrivial if, and only if, there exists $u_n\xrightarrow{n\to\infty}u$ such that $G(u_n)\neq G(u),$ for any $n\in\N.$
\end{defi}

\noindent
When $\ovl u$ is a trivial critical point of $G,$ any concave function $\vphi\in C^0\big([0,r_0)\big)\cap C^1(0,r_0)$ such that $\vphi^\p>0$ and $\vphi(0)=0$ is desingularizing at $\ovl u.$
\medskip \\
An immediate consequence of the KL inequality is a local and strong version of Sard's theorem.
\begin{rmk}[\textbf{Local finiteness of critical values}]
\label{rmkV}
Let $G\in C^1(\R^N;\R)$ and $\ovl u\in \R^N.$ Assume that $G$ satisfies the KL property at $\ovl u$ on $B(\ovl u,\eta).$ Then
\begin{gather*}
u\in B(\ovl u,\eta) \mbox{ and }\nabla G(u)=0 \implies G(u)=G(\ovl u).
\end{gather*}
\end{rmk}

\noindent
The simplest functions we can think of with respect to the behavior of the solutions to \eqref{second1} are given by functions with linear gradients, that is quadratic forms 
\begin{gather*}
G(u)=\frac{1}{2}\langle Au,u\rangle, \:u\in \R^N, \mbox{ where }A\in\Mr_N(\R), \: A^T=A.
\end{gather*}
When $A\neq 0,$ it is easy to establish directly that $\vphi(s)=\sqrt{\frac{1}{|\lambda|}s}$ (where $\lambda$ is a nonzero eigenvalue with smallest absolute value) provides a desingularizing function. In the subsections to come, we show that the best we can hope in general for a desingularizing function $\vphi$ attached to a $C^2$ function $G$ is precisely a quantitative behavior of square-root type. 

\subsection{Lower bounds for desingularizing functions of potentials having a simple critical point structure}

Our first assumption, formally stated below, asserts that points having critical value must be critical points. The assumption is rather strong in general but it will be complemented in the next section by a far more general result for definable functions.
\begin{gather}
\label{A}
\left\{
\begin{array}{l}
\text{Let } \ovl u\in\crit G. \medskip \\
\text{There exists } \eta>0 \text{ such that for any } u\in B(\ovl u,\eta), \smallskip \\
\big(G(u)=G(\ovl u) \implies \;u\in\crit G\big).
\end{array}
\right.
\end{gather}

\begin{exa}
\label{example}
(a) When $N=1$ and $G\in C^1$ is KL then assumption~\eqref{A} holds.\\
\noindent
[{\small If the result does not hold then there exists a sequence $(x_n)_{n\in\N}$ such that $x_n\xrightarrow[]{n\tends\infty}\ovl u$ and
\begin{gather}
\label{dempropG1}
G(x_n)=G(\ovl u),	\\
\label{dempropG2}
G^\p(x_n)\neq0,
\end{gather}
for any $n\in\N.$ Without loss of generality, we may assume that $(x_n)_{n\in\N}$ is monotone, say decreasing. From
\eqref{dempropG1}--\eqref{dempropG2} and Rolle's Theorem, there exists a sequence $(u_n)_{n\in\N}$ such that
$x_{n+1}<u_n<x_n,$ $G^\p(u_n)=0,	$ 
 $G(u_n)\neq G(\ovl u),$ 
for any $n\in\N.$ Thus $G(u_n)$ are critical values distinct from $G(\ovl u)$ such that $G(u_n)\tends G(\ovl u);$ this contradicts the local finiteness of critical values -- see Remark~\ref{rmkV}.}]
\\
(b) Of course, the result in (a) cannot be extended to higher dimensions. Consider for instance 
\begin{gather*}
G:\R^2\tends\R, \quad G(u_1,u_2)=u_1^2-u_2^2,
\end{gather*}
which is obviously KL. One has $\nabla G(u)=0$ if, and only if, $u=0,$ yet $G(t,-t)=0$ for any $t$ in $\R.$
\\
(c) If $G$ is convex, \eqref{A} holds globally, \textit{i.e.}, with $\eta=\infty.$ [{\small This follows directly from the well-known fact that $G(u)=\min G$ if, and only if, $\nabla G(u)=0.$}]
\end{exa}

\begin{lem}[\textbf{Comparing values growth with gradients growth}]
\label{lemG}
\text{} \\
Let $G\-\in C^{1,1}_\loc(\R^N;\R)$ and $\ovl u\in\crit G.$ Assume there exists $\eps>0$ such that 
\begin{gather*}
u\in B(\ovl u,2\eps)\mbox{ and }G(u)=G(\ovl u) \implies u\in\crit G,
\end{gather*}
in other words~\eqref{A} holds $($with $\eta=2\eps).$ Then there exists $c>0$ such that
\begin{gather}
\label{lemG1}
|G(u)-G(\ovl u)|\ge c\|\nabla G(u)\|^2,
\end{gather}
for any $u\in B(\ovl u,\eps).$
\end{lem}

\begin{proof*}
Working if necessary with $\wt G(u)=G(u)-G(\ovl u),$ we may assume, without loss of generality, that $G(\ovl u)=0.$ Let us proceed in two steps. \\
\textbf{Step 1.}
Let $H\in C^{1,1}\big(\ovl B(\ovl u,2\eps);\R\big)$ with $\ovl u\in\crit H$ and assume further that $H\ge0.$ We claim that there exists $c>0$ such that
\begin{gather}
\label{lemGstep1}
\forall u\in B(\ovl u,\eps), \; H(u)\ge c\|\nabla H(u)\|^2.
\end{gather}
Denote by $L_2$ the Lipschitz constant of $\nabla H$ on $\ovl B(\ovl u,2\eps),$ let $L_1=\vmax_{u\in\ovl B(\ovl u,2\eps)}\|\nabla H(u)\|$ and set $L=L_1+L_2.$ Since,
\begin{gather*}
\big(L_1=0 \; \text{ or } \; L_2=0\big) \implies \nabla H_{|B(\ovl u,\eps)}\equiv0 \implies \eqref{lemGstep1},
\end{gather*}
we may assume that $L_2>0$ and $L_1>0.$ Let $u\in B(\ovl u,\eps).$ We have for any $v\in B(0,2\eps),$
\begin{align*}
	& \; H(v)-H(u)=\int_0^1\langle\nabla H\big((1-t)u +tv\big), v-u \rangle\d t \\
   =	& \; \int_0^1\langle\nabla H\big((1-t)u +tv\big)-\nabla H(u), v-u \rangle\d t + \langle\nabla H(u),v-u\rangle,
\end{align*}
so that for any $v\in B(0,2\eps),$
\begin{gather}
\label{estim}
\Big|H(v)-H(u)-\langle\nabla H(u),v-u\rangle\Big|\le\frac{L_2}2\|v-u\|^2.
\end{gather}
Note that $\left\|\left(u-\frac\eps{L}\nabla H(u)\right)-\ovl u\right\|\le\|u-\ovl u\|+\frac\eps{L}\|\nabla H(u)\|<\eps+\eps\frac{L_1}L<2\eps.$ By convexity, we infer that $\left[u,u-\frac\eps{L}\nabla H(u)\right]\subset B(\ovl u,2\eps).$ It follows that $v=u-\frac\eps{L}\nabla H(u)$ is an admissible choice in \eqref{estim}. Without loss of generality, we may assume that $\eps\le1.$ This leads to
\begin{gather*}
0\le H(v)\le H(u) - \frac\eps{2L} \|\nabla H(u)\|^2.
\end{gather*}
Whence the claim. \\
\textbf{Step~2.} Define for any $u\in\ovl B(\ovl u,2\eps),$ $H(u)=|G(u)|.$ Since $\big(G(u)=0\implies\nabla G(u)=0\big),$ we easily deduce that $H\in C^1_\b\big(B(\ovl u,2\eps);\R\big)$ and for any $u\in B(\ovl u,2\eps),$ $\nabla H(u)=\sign\big(G(u)\big)\nabla G(u).$ Denote by $L_2$ the Lipschitz constant of $\nabla G$ on $\ovl B(\ovl u,2\eps).$ We claim that,
\begin{gather}
\label{lemGstep2}
\|\nabla H(u)-\nabla H(v)\|\le L_2\|u-v\|,
\end{gather}
for any $(u,v)\in B(\ovl u,2\eps)\times B(\ovl u,2\eps).$ Let $(u,v)\in B(\ovl u,2\eps)\times B(\ovl u,2\eps).$ Estimate \eqref{lemGstep2} being clear if $G(u)G(v)\ge0,$ we may assume that $G(u)G(v)<0.$ By the Mean Value Theorem and the assumptions on $G,$ it follows that there exists $t\in(0,1)$ such that for $w=(1-t)u+tv,$ $G(w)=0$ and $\nabla G(w)=0.$ We then infer,
\begin{align*}
	&	\; \|\nabla H(u)-\nabla H(v)\|=\|\nabla G(u)+\nabla G(v)\|\le\|\nabla G(u)\|+\|\nabla G(v)\|	\\
   =	&	\; \|\nabla G(u)-\nabla G(w)\|+\|\nabla G(w)-\nabla G(v)\|							\\
  \le	&	\; L_2\|u-w\|+L_2\|w-v\|=L_2\|u-v\|.
\end{align*}
Hence \eqref{lemGstep2}. It follows that $H\in C^{1,1}\big(\ovl B(\ovl u,2\eps);\R\big)$ and $H$ satisfies the assumptions of Step~1. Applying \eqref{lemGstep1} to $H,$ we get \eqref{lemG1}. This concludes the proof.
\medskip
\end{proof*}

\begin{prop}[\textbf{Lower bound for desingularizing functions}]
\label{propG}
Let $G\in C^{1,1}_\loc(\R^N;\R)$ and let $\ovl u$ be a nontrivial critical point, \textit{i.e.} $\ovl u\in\crit G\setminus\inte\crit G.$ Assume that $G$ satisfies the KL property at $\ovl u$ and that assumption~\eqref{A} holds at $\ovl u.$
\\
Then there exists $\beta>0$ such that for any desingularizing function $\vphi$ of $G$ at $\ovl u,$ 
\begin{gather}
\label{grow_ineq}
\vphi^\p(s)\ge\frac\beta{\sqrt s},
\end{gather}
for any small positive $s.$
\end{prop}

\begin{proof*}
 We may assume $G(\ovl u)=0.$ Combining \eqref{Kurd} and \eqref{lemG1}, we deduce that $\vphi^\p(|G(u)|)\ge\frac1{\|\nabla G(u)\|}\ge\frac\beta{\sqrt{|G(u)|}},$ 
for any $u\in B(\ovl u,\eps)$ such that $G(u)\neq G(\ovl u)$ (Remark~\ref{rmkV}). Changing $G$ into $-G$ if necessary, there is no loss of generality to assume that there exists $u_n$ such that $u_n\tends\ovl u$ with $G(u_n)>0$ (recall $\ovl u$ is a nontrivial critical point). Since $G$ is continuous, this implies by a connectedness argument that for some $\rho$ there exists $r>0$ such that $\big|G\big(B(\ovl u,\rho)\big)\big|\supset(0,r).$ Using the parametrization $s\in(0,r)$ we conclude that
$\vphi^\p(s)\ge\frac\beta{\sqrt s},$ for any $s$ sufficiently small.
\medskip
\end{proof*}

\subsection{Lower bounds for desingularizing functions of definable $\bs{C^2}$ functions}

This part makes a strong use of definability arguments (these are recalled in the last section).

\begin{lem}[\textbf{Lower bounds for desingularizing functions of $\bs{C^2}$ definable functions}]
\label{lower}
Let $G:\O\tends\R$ be a $C^2$ definable function on an open subset $\O\ni 0$ of $\R^N.$ We assume that $0$ is a nontrivial critical point$^($\footnote{Equivalently, we assume that there exists $u_n\xrightarrow{n\to\infty}0$ such that $G(u_n)\neq 0.$}$^)$ and that $G(0)=0.$
\\
Since $G$ is definable it has the KL property$^($\footnote{See Theorem~\ref{KL}.}$^)$ that is, there exist $\eta,r_0>0$ and $\vphi:[0,r_0)\tends\R$ as in Definition~$\ref{KL_property}$ such that 
\begin{gather}
\label{KLineq}
\|\nabla \big( \vphi\circ |G|\big)(u)\|\ge 1,
\end{gather}
for any $u$ in $B(0,\eta)$ such that $G(u)\neq 0.$
\\
Then there exists $c>0$ such that
\begin{gather}
\label{low}
\vphi^\p(s)\ge\frac{c}{\sqrt{s}},
\end{gather}
so that $\vphi(s)\ge 2c\sqrt{s},$ for any small $s>0.$
\end{lem}
\medskip

\begin{proof*}
Let us outline the ideas of the proof: after a simple reduction step, we show that the squared norm of a/the smallest gradient on a level line increases at most linearly with the function values. In the second step, we show that this estimate is naturally linked to the increasing rate of $\vphi$ itself and to property~\eqref{low}. Let $\vphi:[0,r_0)\tends\R$ be any desingularizing function of $G$ at $0$ on $B(0,\eta),$ as in Definition~\ref{KL_property}.
\medskip \\
Changing $G$ in $-G$ if necessary, we may assume by Definition~\ref{nonflatcrit}, without loss of generality, that there exists a sequence $(u_n)_n$ such that $u_n\xrightarrow{n\to\infty}0$ and $G(u_n)>0,$ for any $n\in\N.$ Let us proceed with the proof in three steps.
\medskip \\
\textbf{Step 1.}
We first modify the function $G$ as follows. Let $\rho\in C^2(\R^N;[0,1])$ be a semi-algebraic function such that
\begin{gather*}
\left\{
\begin{array}{l}
\supp \rho\subset B(0,\eta)\subset\O, \medskip \\
\rho(x)=1, \text{ if }x\in B\left(0,\frac\eta2\right).
\end{array}
\right.
\end{gather*}
Let us define $\wh G$ on $\R^N$ by
\begin{gather*}
\wh G(u)=
	\begin{cases}
		\rho(u)G(u)+\dist\left(u,B\left(0,\frac\eta2\right)\right)^3,	&	\text{if } u\in\O,				\medskip \\
		0,											&	\text{if } u\in\R^N\setminus\O.
	\end{cases}
\end{gather*}
It follows that $\wh G\in C^2(\R^N;\R),$ leaves the set of desingularizing functions at $0$ unchanged, has compact lower level sets and is definable in the same structure (recall Definition~\ref{domin} (iii)). Finally, we obviously have,
\begin{gather}
\label{nonflat2}
u_n\xrightarrow{n\to\infty}0 \; \text{ with } \; \wh G(u_n)>0, \,\; \forall n\in\N.
\end{gather}
Without loss of generality, we may assume that $\eta\le1$ and $r_0\le\frac{\eta^3}8.$ Let $u\in\R^N\setminus B(0,\eta).$ One has,
\begin{gather*}
\wh G(u)=\dist\left(u,B\left(0,\frac\eta2\right)\right)^3=\left(\|u\|-\frac\eta2\right)^3\ge\frac{\eta^3}8\ge r_0.
\end{gather*}
It follows that,
\begin{gather}
\label{beg}
\inf_{u\in B(0,\eta)\cap[\wh G=r]}\|\nabla\wh G(u)\|=\min_{u\in[\wh G=r]}\|\nabla\wh G(u)\|, \;\; \forall r\in(0,r_0).
\end{gather}
\textbf{Step 2.}
For $r>0,$ we introduce 
\begin{gather*}
(P_r)\qquad \psi(r)=\min\left\{\frac12\|\nabla\wh G(u)\|^2; \;u\in\R^N, \, \wh G(u)=r\right\}.
\end{gather*}
Since the set of critical values of a definable function is finite and since the level sets are compact, we may choose, if necessary, $r_0$ so that $\psi>0$ on $(0,r_0)$ (the fact that $0$ is a nontrivial critical point excludes the case when $\psi$ vanishes around $0).$ If we denote by $S(r)$ the nonempty compact set of solutions to $(P_r),$ one easily sees that
\begin{gather*}
S:(0,r_0)\rightrightarrows\R^N,
\end{gather*}
is a definable point-to-set mapping -- this follows by a straightforward use of quantifier elimination (\textit{i.e.}, by the use of Definition~\ref{domin}). Using the Definable Selection Lemma (Lemma \ref{selec}), one obtains a definable curve $u:(0,r_0)\tends\R^N$ such that $u(r)\in S(r),$ for any $r\in (0,r_0).$ Finally, using the Monotonicity Lemma (Lemma~\ref{mono}) repeatedly on the coordinates $u_i$ of $u,$ one can shrink $r_0$ so that $u$ is actually in $C^1((0,r_0);\R^N).$
\medskip \\
Fix now $r$ in $(0,r_0).$ Since $r$ is noncritical the problem $(P_r)$ is qualified and we can apply Lagrange's Theorem for constrained problems.
This yields the existence of a real multiplier $\lambda(r)$ such that 
\begin{gather}
\label{opt}
\nabla^2\wh G(u(r))\nabla\wh G(u(r))-\lambda(r)\nabla\wh G(u(r))=0,
\end{gather}
with of course $\wh G(u(r))=r.$
\medskip \\
Note that for any $r\in(0,r_0),$ $\nabla\wh G(u(r))\neq0$ (as seen at the beginning of this step) so that $\lambda(r)$ is an actual eigenvalue of $\nabla^2\wh G(u(r)).$ Since $\wh G$ is $C^2,$ the curve $\nabla^2\wh G(u(r))$ is bounded in the space of matrices $\Mr_N(\R).$ Since eigenvalues depend continuously on operators, one deduces from the previous remarks that there exists $\ovl \lambda\ge 0$ such that
\begin{gather*}
|\lambda(r)|\le\ovl\lambda,\:\forall r\in (0,r_0).
\medskip
\end{gather*}
Multiplying \eqref{opt} by $u^\p(r)$ gives $\langle\nabla^2\wh G(u(r))\nabla\wh G(u(r)),u^\p(r)\rangle = \lambda(r)\langle\nabla\wh G(u(r)),u^\p(r)\rangle,$ which is nothing else than 
\begin{gather*}
\frac12\frac{\d}{\d r}\|\nabla\wh G(u(r))\|^2=\lambda(r)\frac{\d}{\d r}\wh G(u(r)).
\end{gather*}
Since $\wh G(u(r))=r,$ one has 
\begin{gather*}
\frac12\frac{\d}{\d r}\|\nabla\wh G(u(r))\|^2=\lambda(r),
\end{gather*}
so after integration on $[s,r]\subset(0,r_0),$ one obtains 
\begin{gather}
\label{inter}
\left|\|\nabla\wh G(u(r))\|^2-\|\nabla\wh G(u(s))\|^2\right|
=2\left|\int_s^r\lambda(\tau)d\tau\right|\le2\ovl\la|r-s|\xrightarrow{r,s\to0}0.
\end{gather}
It follows that $\left(\|\nabla\wh G(u(r))\|^2\right)_{s>0}$ is a Cauchy's family, so that the limit $\ell$ of $\|\nabla\wh G(u(s))\|^2$ as $s$ goes to zero exists in $[0,\infty).$ We recall that by assumption \eqref{nonflat2}, $u_n\xrightarrow{n\to\infty}0,$ $\wh G(u_n)>0$ and $\nabla\wh G(u_n)\xrightarrow{n\to\infty}0.$ Now, setting $r_n=\wh G(u_n),$ one has by definition of $u(r_n),$ $\|\nabla\wh G(u_n)\|\ge\|\nabla\wh G(u(r_n))\|.$ This implies that $\ell=0$ and as a consequence \eqref{inter} yields
\begin{gather}
\frac12\|\nabla\wh G(u(r))\|^2= \int_0^r\lambda(\tau)d\tau\le\ovl \lambda r,
\end{gather}
in other words
\begin{gather}
\label{majpsi}
\psi(r)\le\ovl\lambda r, \;\; \forall r\in (0,r_0).
\end{gather}
\textbf{Step 3.}
Let us now conclude. By KL inequality one has for any $r\in(0,r_0),$
\begin{gather}
\label{ccl}
\vphi^\p(r)\ge\frac{1}{\|\nabla\wh G(u)\|}, \quad \forall u\in B(0,\eta)\cap [G=r].
\end{gather}
As a consequence, we can use \eqref{beg} in \eqref{ccl} and the linear estimate \eqref{majpsi} above to conclude as follows:
\begin{eqnarray*}
\vphi^\p(r)		&	\ge	&	\frac{1}{\inf\left\{\|\nabla\wh G(u)\|; \; u \in B(0,\eta)\cap[\wh G=r]\right\}}	\\
			&	 =	&	\frac{1}{\min\left\{\|\nabla\wh G(u)\|; \; u \in [\wh G=r]\right\}}				\\
			&	\ge	&	\frac{1}{\sqrt{2\psi(r)}}										\\
			&	\ge	&	\frac{c}{\sqrt{r}},
\end{eqnarray*}
for any $r\in(0,r_0),$ with $c=\left(\sqrt{2\ovl \lambda}\right)^{-1}.$ Hence \eqref{low}.
\medskip
\end{proof*}

\begin{rmk}
(a) Note that if $G\not\in C^2$ then~\eqref{low} does not hold. Indeed, take $G(u)=u^\frac32$ and $\vphi(s)=s^\frac23$ as a (semi-algebraic) counter-example. \\
(b) When we omit the assumption that $0$ is a nontrivial critical point, \textit{i.e.} $0\in\inte\crit G,$ then $G$ vanishes in a neighborhood of $0.$ In that case, the result is not true in general since any concave increasing function adequately regular is desingularizing for $G.$ However a function $\vphi(s)=c\sqrt{s}$ can still be chosen as a desingularizing function. 
\\
Hence, {\em for an arbitrary $C^2$ definable function, we can always assume that for any critical point, the corresponding desingularizing function satisfies $\vphi^\p(s)\ge c\frac{1}{\sqrt{s}}$ $($locally for some positive constant $c).$}
\end{rmk}

\section{Damped second order gradient systems}
\label{quasi_gradient}

\subsection{Quasi-gradient structure and KL inequalities}
\label{Lyapunov}

\begin{defi}
\label{defsqg}
Let $\Gamma$ be a nonempty closed subset of $\R^N$ and let $F:\R^N\tends\R^N$ be a locally Lipschitz continuous mapping. 
\begin{itemize}
\item[(i)]
We say that the first order system
\begin{gather}
\label{sqg}
u^\p(t)+F\big(u(t)\big)=0, \; t\in\R_+,
\end{gather}
has a \textit{quasi-gradient structure for} $E$ \textit{on} $\Gamma,$ if there exist a differentiable function $E:\R^N\tends\R$ and
$\alpha_\Gamma=\alpha>0$ such that
\begin{align}
\label{ac}
& \text{\textbf{(angle condition)}} & \big\langle\nabla E(u),F(u)\big\rangle\ge\alpha\,\|\nabla E(u)\|\,\|F(u)\|,	\text{ for any } u\in \Gamma, \\
\label{zero}
&\text{\textbf{(rest-points equivalence)}} & \crit E\, \cap \, \Gamma=F^{-1}(\{0\})\, \cap\, \Gamma.
\end{align}
 \item[(ii)]
 Equivalently a vector field $F$ having the above properties is said to be \textit{quasi-gradient for} $E$ \textit{on}~$\Gamma.$
\end{itemize}
\end{defi}

\noindent
The following result involves classical material and ideas, yet, the fact that an asymptotic alternative can be derived in this setting does not seem to be well-known (see however \cite{MR2674728} in a discrete context).

\begin{thm}[\textbf{Asymptotic alternative for quasi-gradient fields}]
\label{thmsqg}
Let $F:\R^N\tends\R^N$ be a locally Lipschitz mapping that defines a quasi-gradient vector field for $E$ on $\R^N,$ for some differentiable function $E:\R^N\tends\R.$ Assume further that the function $E$ is KL. Let $u$ be any solution to \eqref{sqg}. Then,
\begin{itemize}
\item[$($i\,$)$]
either $\|u(t)\|\xrightarrow{t\to\infty}\infty,$
\item[$($ii\,$)$]
or $u$ converges to a singular point $u_{\infty}$ of $F$ as $t\tends\infty.$
\end{itemize}
When $($ii\,$)$ holds then $u^\p\in L^1\big((0,\infty);\R^N\big)$ and $u^\p(t)\xrightarrow{t\to\infty}0.$ Moreover, we have the following estimate,
\begin{gather}
\label{automaj}
\|u(t)-u_\infty\|\le\frac1\alpha\vphi\big(E(u(t))-E(u_\infty)\big),
\end{gather}
where $\vphi$ is a desingularizing function of $E$ at $u_\infty$ and $\alpha$ is the constant in~\eqref{ac}.
\end{thm}

\begin{proof*} We assume that (i) does not hold, so there exist $u_\infty\in \R^N$ and a sequence $s_n\nearrow\infty$ such that $u(s_n)\xrightarrow{n\tends\infty}u_\infty.$ Note that by continuity of $E,$ one has $E\big(u(s_n)\big)\xrightarrow{n\tends\infty}E(u_\infty).$ Observe also that from equation \eqref{sqg} and the angle condition~\eqref{ac}, one has for any $t\ge0,$
\begin{eqnarray}
\nonumber
\frac{\d}{\d t}\big(E\circ u\big)(t)		&	=	&	\big\langle\nabla E\big(u(t)\big),u^\p(t)\big\rangle		\\
\nonumber
							&	=	&	-\big\langle\nabla E\big(u(t)\big),F\big(u(t)\big)\big\rangle	\\ 
\label{rest}
							&     \le	&	-\alpha \|\nabla E(u(t))\| \, \|F(u(t))\|,
\end{eqnarray}
and thus the mapping $t\longmapsto E(u(t))$ is nonincreasing, which implies
\begin{gather*}
\lim_{t\to\infty} E(u(t))=E(u_\infty).
\end{gather*}
Note that if $E(u(\ovl t))=E(u_\infty)$ for some $\ovl t,$ one would have $\frac{\d}{\d t}\big(E\circ u\big)(t)=0$ for any $t>\ovl t,$ which would in turn imply, by \eqref{rest}, that $\|\nabla E(u(t))\|\,\|F(u(t))\|=0,$ for any such $t.$ In view of the rest point equivalence \eqref{zero}, this would mean that $F(u(t))=0,$ hence by uniqueness of solution curves, that $u(t)=u_\infty$ for any $t\ge0.$ We can thus assume without loss of generality that 
\begin{gather}
 E(u(t))>E(u_\infty), \: \forall t\ge0.
\end{gather}
Let $t_0>0$ be such that $u(t_0)\in B\left(u_\infty,\frac{\eta}2\right)$ and $\vphi\big(E\big(u(t_0)\big)-E(u_\infty)\big)\in \left(0,\frac{\eta\alpha}2\right),$ where $\alpha>0$ is the constant in \eqref{ac} [in view of our preliminary comments and of the continuity of $E$ such a $t_0$ exists]. By continuity of $u,$ there exists $\tau>0$ such that for any $t\in[t_0,t_0+\tau),$ $u(t)\in B(u_\infty,\eta).$ So we may define $T\in(t_0,\infty]$ as
\begin{gather*}
T=\sup\Big\{t>t_0\, ;\, \forall s\in[t_0,t),\ u(s)\in B(u_\infty,\eta)\Big\}.
\end{gather*}
By \eqref{rest}, the Kurdyka-{\L}ojasiewicz inequality \eqref{Kurd} and equation \eqref{sqg}, we have for any $t\in(t_0,T),$
\begin{align}
\nonumber
	&	-\frac{\d}{\d t}\left(\vphi\circ\Big(E\big(u(\:.\:)\big)-E\big(u_\infty\big)\Big)\right)(t)							\\
\nonumber
   = 	&	-\vphi^\p\big(E\big(u(t)\big)-E(u_\infty)\big)\frac{\d}{\d t}\big(E\circ u\big)(t)								\\
\nonumber
 \ge	&	\;\alpha\:\vphi^\p\big(E\big(u(t)\big)-E(u_\infty)\big)\,\big\|\nabla E\big(u(t)\big)\big\|\,\big\|F\big(u(t)\big)\big\|	\\ \nonumber
   =	&	\;\alpha\:\big\|F\big(u(t)\big)\big\|\,\big\|\nabla\big(\vphi\circ \big(E(\:.\:)-E(u_\infty)\big)\big(u(t)\big)\big\| 		\\
\label{l-un}
 \ge	&	\;\alpha\|u^\p(t)\|. 
\end{align}
It follows from the above estimate that
\begin{gather}
\label{proofthmsqg2}
\|u(t)-u(t_0)\|\le\vint_{t_0}^t\|u^\p(s)\|\d s\le\frac{\vphi\big(E\big(u(t_0)\big)-E(u_\infty)\big)}\alpha<\frac{\eta}2,
\end{gather}
for any $t\in(t_0,T).$ We claim that $T=\infty.$ Indeed, otherwise $T<\infty$ and \eqref{proofthmsqg2} applies with $t=T.$ Hence,
\begin{gather*}
\|u(T)-u_\infty\|\le\|u(T)-u(t_0)\|+\|u(t_0)-u_\infty\|<\eta.
\end{gather*}
Then $u(T)\in B(u_\infty,\eta),$ which contradicts the definition of $T.$ As a consequence the curve $u^\p$ belongs to $L^1\big((t_0,\infty);\R^N\big)$ by \eqref{proofthmsqg2} and the curve $u$ converges to $u_\infty$ by Cauchy's criterion. Finally since $0$ must be a cluster point of $u^\p$ \big(recall indeed $\int_0^\infty\|u^\p(t)\|\d t<\infty$ and $u^\p$ is uniformly continuous by \eqref{sqg}\big), one must have $F(u_\infty)=0.$ The announced estimate follows readily from \eqref{proofthmsqg2} and the fact that $T=\infty.$
\medskip
\end{proof*}

\begin{cor}
\label{corsqg}
Let $F:\R^N\tends\R^N$ be locally Lipschitz continuous and assume that for any $R>0$ the mapping $F$ defines a quasi-gradient vector field for some differentiable function $E_R:\R^N\tends\R$ on $\ovl B(0,R).$ Assume further that each of the functions $E_R$ is KL. \\
Let $u$ be any bounded solution to \eqref{sqg}. Then $u$ converges to a singular point $u_\infty$ of $F,$ $u^\p$ is integrable and converges to $0.$ In particular, if we take $R\ge\sup \big\{\|u(t)\|; \; t\in [0,\infty)\big\},$ we have the following estimate,
\begin{gather}
\|u(t)-u_\infty\|\le\frac1{\alpha_R}\vphi\Big(E_R(u(t))-E_R(u_\infty)\Big),
\end{gather}
where $\vphi$ is a desingularizing function of $E_R$ at $u_\infty$ and $\alpha_R$ is the constant in~\eqref{ac}, for the ball $\ovl B(0,R).$
\end{cor}

\begin{proof*}
Take $R\ge\sup \big\{\|u(t)\|; \,t\in [0,\infty)\big\}$ and observe that the previous proof may be reproduced as it is: just replace $E$ by $E_R.$
\medskip
\end{proof*}

\subsection{Convergence rate of quasi-gradient systems and worst-case dynamics}

To simplify our presentation we consider first a proper gradient system:
\begin{gather}
\label{gradsystem}
u^\p(t)+\nabla E(u(t))=0,
\end{gather}
where $E:\R^N\tends\R$ is a twice continuously differentiable KL function. We assume that $u$ is bounded so, by virtue of our previous considerations, the curve converges to some critical point $u_\infty$ of $E.$ Observe that if $u_\infty$ is a trivial critical point, one actually has $u(0)=u_\infty$ and the asymptotic study is trivial.
\medskip \\
We thus assume $u_\infty$ to be nontrivial, and we denote by $\vphi$ a desingularizing function of $E$ at $u_\infty.$ We set
\begin{gather*}
\psi=\vphi^{-1},
\end{gather*}
whose domain is denoted by $[0,a),$ (with $a\in(0,\infty])$ and we consider the {\em one-dimensional worst-case gradient dynamics} (see \cite{HDR}):
\begin{gather}
\label{worstcase}
\nu^\p(t)+\psi^\p(\nu(t))=0,\quad \nu(0)=\nu_0\in (0,a).
\end{gather}
We shall assume that
\begin{align}
\label{nondeg}
\vphi^\p(s)\ge\frac{c}{\sqrt{s}}, \text{ on }(0,r_0), &
\end{align}
which implies that solutions $\nu$ to \eqref{worstcase} are globally defined on $[0,\infty)$ and satisfy $\vlim_{t\nearrow\infty}\nu(t)=0$ with $\nu(t)\ge\nu_0e^{-c_0t},$ for any $t\ge0$ (and for some $c_0>0).$ Uniqueness holds by concavity of $\vphi.$ Finally, note that if $E$ is a $C^2$ definable function then $\vphi$ can be chosen to be $C^2,$ strictly concave and satisfying \eqref{nondeg} (Remark~\ref{rmkKL} (c) and Lemma~\ref{lower}).
\medskip \\
\textbf{Radial functions and worst-case dynamics.} A full justification of the terminology ``worst-case dynamics" is to be given further, but at this stage one can observe that $E$ could be taken of the form
\begin{gather*}
E_{\mbox{\small rad}}(u)=\vphi^{-1}(\|u-u_\infty\|), \text{ with } u\in B(u_\infty,\eta) \;(\eta>0),
\end{gather*}
provided that $\vphi^{-1}$ is smooth enough. In that case $\vphi$ is clearly desingularizing and the solutions of the gradient system \eqref{gradsystem} are radial in the sense that they are of the form$^($\footnote{Just use the formula in \eqref{gradsystem}.}$^)$
\begin{gather}
\label{rad}
u(t)=u_\infty+\nu(t)\frac{u_0-u_\infty}{\|u_0-u_\infty\|},
\end{gather}
where $\nu$ is a solution to \eqref{worstcase}. In this case, the dynamics \eqref{worstcase} exactly measures the convergence rates for \eqref{gradsystem}, since one has for any $t\ge 0$ and any $u_0$ such that $\nu(0)=\|u_0-u_\infty\|,$
\begin{align}
& E_{\mbox{\small rad}}(u(t))=\psi(\nu(t)), \\
& \|u(t)-u_\infty\|=\nu(t).
\end{align} 
We are about to see that this behavior in terms of convergence rate is actually the worst we can expect. 

\begin{rmk}
\label{vitesse}
(a) As can be seen below, the worst-case gradient system is introduced to measure the rate of convergence of solutions for large $t.$ Since nontrivial solutions to \eqref{worstcase} have the same asymptotic behavior (they are, indeed, all of the form $\nu_1(t)= \nu(t+t_0)$ where $t_0$ is some real number), the choice of the initial condition $\nu(0)$ in $(0,a)$ can be made arbitrarily. \\
(b) The above rewrites $\nu^\p(t)\vphi^\p\big(\vphi^{-1}(\nu(t))\big)=-1.$ Thus if $\mu$ denotes an antiderivative of 
$\vphi^\p\circ \vphi^{-1},$ one has $\nu(t)=\mu^{-1}(-t+a_0)$ (where $a_0$ is a constant), for any $t>0$ large enough. \\
(c) In general, the explicit integration of such a system depends on the integrability properties of $\psi$ and on the fact that $\vphi^\p\circ \vphi^{-1}$ admits an antiderivative in a closed form. \\
For instance if $\vphi(s)=(\frac{s}{c})^{\theta},$ with $c>0$ and $\theta\in\left(0,\frac12\right),$ then $\psi(s)=cs^{\frac{1}{\theta}}$ and
\begin{gather*}
\nu^\p(t)+\frac{c}{\theta}\,\nu(t)^{\frac{1-\theta}{\theta}}=0,\quad \nu(0)\in (0,a).
\end{gather*}
Thus by integration
\begin{gather*}
\frac{\d}{\d t}\nu^{1-\frac{1-\theta}{\theta}}(t)=\frac{\d}{\d t}\nu^{\frac{-1+2\theta}{\theta}}(t)=c_1,
\end{gather*}
with $c_1>0.$ As a consequence,
\begin{gather*}
\nu(t)=\big(c_2+c_1t\big)^{-\frac{\theta}{1-2\theta}},
\end{gather*}
with $c_2>0.$ When $\theta=\frac{1}{2}$ one easily sees that $\nu(t)=\nu(0)\exp\left(-2ct\right).$
\end{rmk}

\begin{thm}[\textbf{The worst-case rate and worst-case one-dimensional gradient dynamics}]
\label{worst}
\text{} \\
Let $E\in C^2(\R^N;\R)$ be a KL function, let $u$ be a bounded solution to \eqref{gradsystem} and let $u_\infty\in\crit E$ satisfying $u(t)\xrightarrow{t\to\infty}u_\infty$ $($such a $u_\infty$ exists by Theorem~$\ref{thmsqg}).$ Then for any $t$ large enough,
\begin{gather}
\label{majval}
E(u(t))-E(u_\infty)\le\psi(\nu(t)),
\end{gather}
and
\begin{gather}
\label{majgrad}
\|u(t)-u_{\infty}\|\le\nu(t),
\end{gather}
where $\nu$ is a solution to \eqref{worstcase}.
\end{thm}

\begin{proof*}
Without loss of generality, we may assume that $E(u_\infty)=0.$ From the previous results, we know that for any $t\ge t_0,$ we have $u(t)\in B(u_\infty,\eta)$ and $E(u(t))\in(0,r_0),$ so that the KL inequality gives (see Theorem~\ref{thmsqg} and \eqref{l-un}):
\begin{gather*}
-\frac{\d}{\d t}\big(\vphi\circ E(u)\big)(t)\ge\|u^\p(t)\|.
\end{gather*}
Set $z(t)=E(u(t)).$ Since $\frac{\d}{\d t}(E\circ u)(t)=-\|u^\p(t)\|^2,$ one has $-\frac{\d}{\d t}(\vphi\circ z)(t)\ge\sqrt{-z^\p(t)},$ or equivalently
\begin{gather*}
\vphi^\p\big(z(t)\big)^2z^\p(t)\le-1.
\end{gather*}
Consider now the worst-case gradient system with initial condition $\nu(t_0)=\vphi\big(E(u(t_0))\big)$ and set
$z_a(t)=\psi(\nu(t))=\vphi^{-1}(\nu(t)),$ for $t\ge t_0.$ The system \eqref{worstcase} becomes
$\vphi^\p(z_a(t))z^\p_a(t)+\frac{1}{\vphi^\p(z_a(t))}=0,$ \textit{i.e.}, $\vphi^\p(z_a(t))^2z^\p_a(t)=-1.$ 
If $\mu$ is an antiderivative of ${\vphi^\p}^2$ on $(0,r_0),$ it is an increasing function and one has 
\begin{gather*}
\frac{\d}{\d t}(\mu\circ z)(t)=\vphi^\p\big(z(t)\big)^2z^\p(t)\le-1=\vphi^\p\big(z_a(t)\big)^2z_a^\p(t)=\frac{\d}{\d t}(\mu\circ z_a)(t),
\end{gather*}
and $\mu(z(t_0))=\mu(z_a(t_0)).$
As a consequence, $\mu(z(t))\le \mu(z_a(t)),$ hence $z(t)\le z_a(t)$ for any $t\ge t_0,$ which is exactly \eqref{majval}. Using \eqref{automaj}, we conclude by observing that 
\begin{gather*}
\|u(t)-u_\infty\|\le\vphi(E(u(t)))\le\vphi(z_a(t))=\nu(t).
\end{gather*}
The theorem is proved.
\medskip
\end{proof*}

\begin{rmk}
Observe that in the case of a desingularizing function of power type (see Remark~\ref{vitesse} (c)), we recover well-known estimates \cite{MR1829143}.
\end{rmk}

\begin{thm}[\textbf{The worst-case one-dimensional gradient dynamics for quasi-gradient systems}]
\label{worstquas}
\text{} \\
Let $F:\R^N\tends\R^N$ be a locally Lipschitz continuous mapping that defines a quasi-gradient vector field for some function $E\in C^2(\R^N;\R)$ on $\ovl B(0,R),$ for any $R>0.$ Assume further that the function $E$ is KL and that for any $R>0,$ there exists a positive constant $b>0$ such that 
\begin{gather}
\label{asfast}
\|\nabla E(u)\|\le b \|F(u)\|,
\end{gather}
for any $u\in\ovl B(0,R).$ Assume further that for a given initial data $u_0\in\R^N$ the solution $u$ to \eqref{sqg} converges to some rest point~$u_\infty.$ Denote by $\vphi$ some desingularizing function for $E$ at~$u_\infty.$
\\
Then there exist some constants $c,d>0, t_0\in \R$ such that
\begin{gather}
\label{majgrad1}
\|u(t)-u_{\infty}\|\le d\nu\left(ct+t_0\right),
\end{gather}
where $\nu$ is a solution to \eqref{worstcase}.
\end{thm}

\begin{proof*}
Combining the techniques used in Theorems~\ref{thmsqg} and~\ref{worst}, the proof is almost identical to that of Theorem~\ref{worst}. Without loss of generality, we may assume that $E(u_\infty)=0.$ We simply need to check the following inequality which is itself a consequence of the assumption \eqref{asfast} applied with $R=\vsup_{t>0}\|u(t)\|.$
\begin{eqnarray*}
-\frac{\d}{\d t}(E\circ u)(t)	&	=	&	-\langle u^\p(t),\nabla E(u(t))\rangle	\\
					&	\le	&	\|F(u(t))\| \, \|\nabla E(u(t))\|		\\
					&	\le	&	b \, \|F(u(t))\|^2					\\
					&	\le	&	b \, \|u^\p(t)\|^2.
\end{eqnarray*}
From \eqref{l-un} one has $-\frac{\d}{\d t}(\vphi\circ E)(u(t)) \ge \alpha \|u^\p(t)\|,$ for any $t$ sufficiently large. Setting $z(t)=E(u(t)),$ one obtains $-\frac{\d}{\d t}(\vphi\circ z)(t)\ge\frac\alpha{\sqrt b}\sqrt{-z^\p(t)}.$ The conclusion follows as before by using a reparametrization of \eqref{worstcase}.
\medskip
\end{proof*}
 
\begin{rmk}
Assumption \eqref{asfast} is of course necessary and simply means that the vector field $F$ drives solutions to their rest points at least ``as fast as $\nabla E$" (see also \cite{MR2572850}).
\end{rmk}

\subsection{Damped second order systems are quasi-gradient systems}
\label{S:desingularizing}

As announced earlier our approach to the asymptotic behavior of damped second order gradient system is based on the observation that \eqref{second1} can be written as a system having a quasi-gradient structure. For $G\in C^2(\R^N;\R),$ let us define $\vF:\R^N\times\R^N\tends\R^N$ by
\begin{gather*}
\vF (u,v) = \big(-v,\,\gamma v + \nabla G(u) \big).
\end{gather*}
Then \eqref{second1} is equivalent to 
\begin{gather}
\label{FOS}
U^\p(t)+\vF\big(U(t)\big)=0, \quad t\in\R_+, \text{ with } U=(u,v).
\end{gather}
As explained in the introduction the total energy function $E_T(u,v)=G(u)+\frac12\|v\|^2$ (sum of the potential energy and the kinetic energy) is a Liapunov function for our dynamical system \eqref{second1}. Formally
\begin{gather*}
\langle \nabla E_T(u,v) , \vF (u,v) \big\rangle = \gamma\|v\|^2.
\end{gather*}
From the above we see, that the damped system \eqref{second1} is not quasi-gradient for $E_T$ since one obviously has a degeneracy phenomenon
\begin{gather}
\label{degen}
\big\langle \nabla E_T(u,v) , \vF (u,v) \big\rangle=0 \; \text{ whenever } \; v=0,
\end{gather}
where in general $\nabla E_T(u,v)\neq 0$ and $\vF (u,v)\neq0.$ 
\\
The idea that follows consists in continuously deforming the level sets of $E_T,$ through a family of functions:
\begin{gather*}
\vE_\la:\R^N\times\R^N\tends \R \; \text{ with } \; \vE_0=E_T \;\; (\lambda \text{ denotes here a positive parameter)},
\end{gather*}
so that the angle formed between each of the gradients of the resulting functions $\vE_\la, \la>0$ and the vector $\vF$ remains far away from $\pi/2.$ In other words we seek for functions making $\vF$ a quasi-gradient vector field. 

\begin{prop}[\textbf{The second order gradient systems are quasi-gradient systems}]
\label{P:lyapunov}
Let $G\in C^2(\R^N;\R)$ and let $\gamma>0.$ For $\la >0,$ define $\vE_\la\in C^1(\R^N\times\R^N;\R)$ by
\begin{gather*}
\vE_\la(u,v) = \left(\frac12 \| v\|^2 + G(u) \right) + \la \langle\nabla G(u), v\rangle.
\end{gather*}
For any $R>0,$ there exists $\la_0>0$ satisfying the following property. For any $\la\in (0,\la_0],$ there exists $\alpha>0$ such that
\begin{gather}
\label{cond1}
\big\langle \nabla\vE_\la(u,v) , \vF (u,v) \big\rangle \ge \alpha \, \| \nabla\vE_\la(u,v) \| \, \| \vF(u,v)\|,
\end{gather}
for any $(u,v)\in\ovl B(0,R)\times\R^N.$ Furthermore,
\begin{gather}
\label{cond2}
\crit\vE_\la\cap\left(\ovl B(0,R)\times\R^N\right)=\vF^{-1}(\{0\})\cap\left(\ovl B(0,R)\times\R^N\right),
\end{gather}
for any $\la\in[0,\la_0].$
\end{prop}

\begin{proof*}
For each $(u,v)\in \R^N\times\R^N,$ we have
$\nabla \vE_\la(u,v) = \left(\nabla G(u) + \la \, \nabla^2 G(u) v,\, v + \la\nabla G(u) \right).$ Let $R>0$ be given and let $M=\max\big\{ \|\nabla^2 G(u)\|; \; u\in\ovl B(0,R)\big\}.$ Choose $\la_0>0$ small enough to have
\begin{gather*}
\gamma-\left(M+\frac{\gamma^2}2\right)\la_0>0.
\end{gather*}
Let $\la\in(0,\la_0].$ Then for any $(u,v)\in\ovl B(0,R)\times\R^N,$ we obtain by Young's inequality,
\begin{align}
\nonumber
\big\langle\nabla \vE_\la(u,v),\vF (u,v)\big\rangle
& =\gamma \|v\|^2-\la\,\langle\nabla^2G(u)v,v\rangle+\la\,\langle\nabla G(u),\gamma v\rangle+\la\,\|\nabla G(u)\|^2 \\
\nonumber
& \ge\left(\gamma-M\la_0-\frac{\la_0}{2}\gamma^2\right)\,\|v\|^2+\frac{\la}{2}\,\|\nabla G(u)\|^2 \\
\label{crit}
& \ge \alpha_0\,(\|v\|^2+\|\nabla G(u)\|^2),
\end{align}
where $\alpha_0=\min\left\{\gamma-\left(M+\frac{\gamma^2}2\right)\la_0,\,\frac\la2\right\}>0.$ Moreover,
\begin{gather}
\label{proofcon1}
\|\nabla \vE_\la(u,v)\|\, \|\vF (u,v) \| \le \frac12\|\nabla \vE_\la(u,v)\|^2+\frac12\|\vF(u,v)\|^2\le C(\|v\|^2+\|\nabla G(u)\|^2 ).
\end{gather}
Combining \eqref{proofcon1} with \eqref{crit}, we deduce that the angle condition \eqref{cond1} is satisfied with $\alpha=\frac{\alpha_0}{C}.$ Finally, the rest point equivalence \eqref{cond2} follows from \eqref{crit}.
\medskip
\end{proof*}

\begin{rmk}
Note that for $\la=0,$ we recover the \textit{total energy} $E_T(u,v)=\vE_0(u,v)=\frac12 \| v\|^2 + G(u).$
\end{rmk}

\medskip

\noindent
The following result is of primary importance: roughly speaking it shows that functions which desingularize the potential $G$ at some critical point $\ovl u,$ also desingularize the energy function $E_T$ and more generally the family of deformed functions $\vE_\la$ at the corresponding critical point
$(\ovl u,0).$ This result implies in turn that the decay rate of the energy is essentially conditioned by the geometry of $G$ as one might expect from a mechanical or an intuitive perspective.
\medskip \\
In the proposition below one needs the kinetic energy to be desingularized by $\vphi.$ This explains our main assumption.

\begin{prop}[\textbf{Desingularizing functions of the energy}]
\label{G_Energy}
Let $G\in C^2(\R^N;\R),$ $\ovl u\in\crit G$ and assume that there exists a desingularizing function $\vphi\in C^1\big((0,r_0);\R_+\big)$ of $G$ at $\ovl u$ on $B(\ovl u,\eta)$ such that $\vphi^\p(s)\ge\frac{c}{\sqrt{s}},$ for any $s\in(0,r_0).$

\noindent
Then there exist $\la_1>0,$ $\eta_1>0$ and $c>0$ such that
\begin{gather}
\label{KL-energy}
\left\|\nabla \left(\vphi \circ \frac12\,|\vE_\la(\:.\:,\:.\:)-\vE_\la(\ovl u,0)|\right)(u,v)\right\|\ge c,
\end{gather}
for any $\la\in [0,\la_1]$ and any $(u,v)\in B(\ovl u,\eta_1)\times B(0,\eta_1)$ such that $\vE_\la(u,v)\neq\vE_\la(\ovl u,0).$
\end{prop}

\begin{proof*}
By standard translation arguments, we may assume without loss of generality that $G(\ovl u)=0$ and $\ovl u=0.$ Then $\vE_\la(0,0)=0$ and \eqref{KL-energy} consists in showing that for some constant $c>0,$
\begin{gather*}
\vphi^\p\left(\frac12\,|\vE_\la(u,v)|\right)\ge \frac{c}{\|\nabla \vE_\la(u,v)\|},
\end{gather*}
for any $\la\in[0,\la_1]$ and any $(u,v)\in B(0,\eta_1)\times B(0,\eta_1)$ such that $\vE_\la(u,v)\neq0.$ Recall that $0\in\crit G.$ Let $M=\max\Big\{\|\nabla^2 G(u)\|; \;u\in\ovl B(0,\eta)\Big\}$ and define
$\la_1=\min\left\{\frac14\,,\frac1{2(M^2+1)}\right\}.$ We have, 
\begin{align}
\label{MAJ}
\|\nabla\vE_\la(u,v)\|^2	&	= \|\nabla G(u)+\la \nabla^2G(u)v\|^2+\|v+\la \nabla G(u)\|^2			\nonumber \\
					&	\ge \|\nabla G(u)\|^2 + \|v\|^2-\la_1(M^2+1)\|v\|^2-2\la_1\|\nabla G(u)\|^2	\nonumber \\
					&	\ge \frac 12\big(\|v\|^2+\|\nabla G(u)\|^2\big),
\end{align}
and in particular,
\begin{gather}
\label{MAJ1}
\|\nabla G(u)\|\le2\|\nabla\vE_\la(u,v)\|,
\end{gather}
for any $\la\in [0,\la_1]$ and any $(u,v)\in\R^N\times\R^N.$ Let now $(\la,u,v)\in[0,\la_1]\times B(0,\eta)\times\R^N$ be such that
$\vE_\la(u,v)\neq0.$ Since $\vphi^\p$ is nonincreasing, we have
\begin{align}
\label{ineg_tri}
\vphi^\p\left(\frac12\,|\vE_\la(u,v)|\right)	&	\ge\vphi^\p\left(\frac12\,|\vE_\la(u,v)-\vE_\la(u,0)|+\frac12\,|\vE_\la(u,0)|\right) \nonumber	\\
								&	\ge\vphi^\p\left(\max\big\{|\vE_\la(u,v)-\vE_\la(u,0)|,|\vE_\la(u,0)|\big\}\right).
\end{align}
Let us first find a lower bound on $\vphi^\p(|\vE_\la(u,0)|).$ Observe that necessarily $\vE_\la(u,0)=G(u)\neq0.$ In particular, $\nabla G(u)\neq0$
(Remark~\ref{rmkV}). We then have by \eqref{Kurd} and \eqref{MAJ1}, $\nabla\vE_\la(u,v)\neq0$ and
\begin{gather}
\label{reg2}
\vphi^\p\big(|\vE_\la(u,0)|\big) = \vphi^\p\big(|G(u)|\big) \ge \frac{1}{\|\nabla G(u)\|} \ge \frac{1}{2\|\nabla \vE_\la(u,v)\|},
\end{gather}
for any $\la\in [0,\la_1]$ and any $(u,v)\in B(0,\eta)\times\R^N$ such that $\vE_\la(u,0)\neq0.$ 
\medskip \\
Let us now estimate $\vphi^\p(|\vE_\la(u,v)-\vE_\la(u,0)|)$ in \eqref{ineg_tri} under the assumption $\vE_\la(u,v)\neq\vE_\la(u,0).$ Cauchy-Schwarz' inequality implies that for any
$\la\in[0,\la_1],$
\begin{gather}
\label{cauchy-schwarz}
|\vE_\la(u,v)-\vE_\la(u,0)|\le \frac12 \big(\|v\|^2+\la_1\|v\|^2+\la_1\|\nabla G(u)\|^2\big).
\end{gather}
Combining \eqref{cauchy-schwarz} with \eqref{MAJ}, we deduce that for any $\la\in [0,\la_1]$ and any $(u,v)\in\R^N\times\R^N,$
\begin{gather}
\label{ligne_v}
|\vE_\la(u,v)-\vE_\la(u,0)|\le(1+\la_1)\|\nabla \vE_\la(u,v)\|^2.
\end{gather}
By continuity of $\nabla G,$ there exists $\eta_1\in(0,\eta)$ such that
\begin{gather*}
\sup\Big\{(1+\la_1)\|\nabla\vE_\la(u,v)\|^2;\, (\la,u,v) \in [0,\la_1]\times B(0,\eta_1)\times B(0,\eta_1)\Big\}<r_0.
\end{gather*}
Using successively the fact that $\vphi^\p$ is nonincreasing and $\vphi^\p(s)\ge\frac{c}{\sqrt{s}},$ it follows from \eqref{ligne_v} that if $(u,v)\in B(0,\eta_1)\times B(0,\eta_1)$ with $\vE_\la(u,v)\neq\vE_\la(u,0)$ then $\nabla \vE_\la(u,v)\neq0$ and
\begin{gather}
\label{reg}
\vphi^\p\big(|\vE_\la(u,v)-\vE_\la(u,0)|\big)\ge\vphi^\p\big((1+\la_1)\|\nabla\vE_\la(u,v)\|^2\big)\ge \frac{c_1}{\|\nabla \vE_\la(u,v)\|},
\end{gather}
where $c_1>0$ is a constant. Finally, inequalities \eqref{reg2} and \eqref{reg} together with \eqref{ineg_tri} yield the existence of a constant $c>0$ such that for any $\la\in [0,\la_1]$ and any $(u,v)\in B(0,\eta_1)\times B(0,\eta_1)$ such that $\vE_\la(u,v)\neq0,$ there holds
$\nabla\vE_\la(u,v)\neq0$ and $\vphi^\p\left(\frac12|\vE_\la(u,v)|\right)\|\nabla\vE_\la(u,v)\|\ge c,$
which is the desired result.
\medskip
\end{proof*}

\section{Convergence results}
\label{convergence}

\noindent
Before providing our last results, we would like to recall to the reader that a bounded trajectory of \eqref{second1} may not converge to a single critical point; finite-dimensional counterexamples for $N=2$ are provided
in~\cite{MR1753136,MR2018329}, in each case the trajectory of \eqref{second1} ends up circling indefinitely around a disk.
\medskip \\
We now proceed to establish a central result whose specialization to various settings will provide us with several extensions of Haraux-Jendoubi's initial work \cite{MR1616968}.

\begin{thm}
\label{thmmain}
Let $G\in C^{2}(\R^N;\R)$ and $(u_0,u_0^\p)\in\R^N\times\R^N$ be a set of initial conditions for \eqref{second1}. Denote by $u\in C^2\big([0,\infty);\R^N)$ the unique regular solution to~\eqref{second1} with initial data $(u_0,u_0^\p).$ Assume that the following holds.
\begin{enumerate}
\item[]
\begin{enumerate}[$1.$]
\item
\textbf{$\bs($The trajectory is bounded$\bs)$} $\vsup_{t>0}\|u(t)\|<\infty.$
\item
\textbf{$\bs($Convergence to a critical point$\bs)$}
$G$ is a KL function. Each desingularizing function $\vphi$ of $G$ satisfies
\begin{gather}
\label{hypo}
\vphi^\p(s)\ge\frac\beta{\sqrt s},
\end{gather}
for any $s\in(0,\eta_0),$ where $\beta$ and $\eta_0$ are positive constants $($see Definition~$\ref{KL_property}).$
\end{enumerate}
\end{enumerate}
Then,
\begin{enumerate}[$(i)$]
\item
$u^\p$ and $u^{\p\p}$ belong to $L^1\big((0,\infty);\R^N\big)$ and in particular $u$ converges to a single limit $u_\infty$ in $\crit G.$
\item
When $u$ converges to $u_\infty,$ we denote by $\vphi$ the desingularizing function of $G$ at $u_\infty.$ One has the following estimate
\begin{gather*}
\|u(t)-u_\infty\|\le c\nu(t),
\end{gather*}
where $\nu$ is the solution of the worst-case gradient system 
\begin{gather*}
\nu^\p(t)+(\vphi^{-1})^\p(\nu(t))=0, \; \nu(0)>0.
\end{gather*}
\end{enumerate}
\end{thm}

\begin{vproof}{of Theorem~\ref{thmmain}.}
Let $G\in C^2(\R^N;\R),$ let $(u_0,u_0^\p)\in\R^N\times\R^N,$ let $u\in C^2\big([0,\infty);\R^N)$ and let $\ovl u\in\R^N.$ Set $U(t)=\big(u(t),u^\p(t)\big),$ $U_0=(u_0,u_0^\p)$ and $\ovl U=(\ovl u,0).$ Let $\vF$ and let $\vE_\la$ be defined as in Subsection~\ref{Lyapunov} and Proposition~\ref{P:lyapunov}, respectively. Note that if $\ovl u\not\in\crit G$ then $\ovl U\not\in\crit\vE_\la$ and $\vphi(t)=ct$ desingularizes $\vE_\la$ at $\ovl U,$ for any $\la\ge0$ \big(Remark~\ref{rmkKL}~(a) and \eqref{cond2}\big). Otherwise, $\ovl u\in\crit G$ and we shall apply Proposition~\ref{G_Energy}. Since $\sup_{t>0}\|u(t)\|<\infty,$ $u^\pp(t)+\gamma u^\p(t)=A(t)$ where $A$ is bounded. Thus, $u^\p(t)=u^\p(0)e^{-\gamma t}+\int_0^t \exp(-\gamma(t-s))A(s)\d s,$ and by a straightforward calculation, $\sup_{t>0}\|u^\p(t)\|<\infty.$ It follows that $\sup_{t>0}\|U(t)\|<\infty.$ Let $R=\sup_{t>0}\|U(t)\|.$ Let $\la_0>0$ and $0<\la_1<\la_0$ be given by Propositions~\ref{P:lyapunov} and \ref{G_Energy}, respectively. Let us fix $0<\la_\star<\la_1$ and let $\alpha>0$ be given by Proposition~\ref{P:lyapunov} for such $\vE_{\la_\star}$ and $R.$ By Proposition~\ref{P:lyapunov}, the first order system
\begin{gather}
\label{proof}
U^\p(t)+\vF\big(U(t)\big)=0, \quad t\in\R_+,
\end{gather}
has a quasi-gradient structure for $\vE_{\la_\star}$ on $\ovl B(0,R)$ (Definition~\ref{defsqg}). Finally, since $G$ has the KL property at $\ovl u,$ $\vE_{\la_\star}$ also has the KL property at $\ovl U$ (Proposition~\ref{G_Energy}). It follows that Theorem~\ref{thmsqg} applies to $U,$ from which $(i)$ follows.
\medskip

\noindent
The estimate part of the proof of $(ii)$ will follow from Theorem~\ref{worstquas}, if we establish that for any $R>0,$ there exists $b>0$ such that for any $(u,v)\in\ovl B(0,R)\times\ovl B(0,R),$
\begin{gather*}
\|\nabla \vE_{\la_\star}(u,v)\|\le b\|\vF(u,v)\|.
\end{gather*}
First we observe that for each $R>0$ and for any $(u,v)\in\ovl B(0,R)\times\ovl B(0,R),$ there exists $k_1\ge 0$ such that 
\begin{gather}
\label{1}
\|\nabla \vE_{\la_\star}(u,v)\|^2 \le k_1\big(\|\nabla G(u)\|^2+\|v\|^2\big).
\end{gather}
This follows trivially by Cauchy-Schwarz' inequality and the fact that $\nabla^2G$ is continuous hence bounded on bounded sets. Fix $\sigma>0$ and recall the inequality $2ab\le\sigma^2a^2+\frac{b^2}{\sigma^2}$ for all real numbers $a,b.$ By Cauchy-Schwarz' inequality and the previous inequality 
\begin{eqnarray*}
\|\vF(u,v)\|^2
	&	=	&	\|v\|^2 +\|\gamma v+ \nabla G(u)\|^2													\\
	&    \ge	&	(1+\gamma^2)\|v\|^2 + \|\nabla G(u)\|^2-2\|\gamma v\|\|\nabla G(u)\|							\\
	&    \ge	&	(1+\gamma^2)\|v\|^2 + \|\nabla G(u)\|^2-\sigma^2\|\gamma v\|^2-\frac{1}{\sigma^2}\|\nabla G(u)\|^2	\\
	&	=	&	(1-(\sigma^2-1)\gamma^2)\|v\|^2 + \left(1-\frac{1}{\sigma^2}\right)\|\nabla G(u)\|^2.
\end{eqnarray*}
Choosing $\sigma>1$ so that $1-(\sigma^2-1)\gamma^2>0$ yields $k_2>0$ such that
$\|\vF(u,v)\|^2 \ge k_2\big(\|\nabla G(u)\|^2+\|v\|^2\big),$ for any $u,v$ in $\R^N.$ Combining this last inequality with \eqref{1}, we obtain $\|\nabla\vE_{\la_\star}(u,v)\|^2\le\frac{k_1}{k_2}\|\vF(u,v)\|^2,$ for any $(u,v)\in\ovl B(0,R)\times\ovl B(0,R).$ Hence the result.
\medskip
\end{vproof}

\begin{rmk} (a) As announced previously convergence rates depend directly on the geometry of $G$ through $\vphi.$\\
(b) The fact that the length of the velocity curve $u^\p$ is finite suggests that highly oscillatory phenomena are unlikely.
\end{rmk}

\section{Consequences}
\label{consequences}

\noindent
In the following corollaries, the mapping $\R_+\ni t\longmapsto u(t)$ is a solution curve of \eqref{second1}.

\begin{cor}[\textbf{Convergence theorem for real-analytic functions \cite{MR1616968}}]
\label{cor1}
Assume that $G:\R^N\tends\R$ is real-analytic and let $u$ be a bounded solution to~\eqref{second1}. Then we have the following result.
\begin{itemize}
\item[$(i)$]
$(u,u^\p)$ has a finite length. In particular $u$ converges to a critical point $u_\infty.$
\item[$(ii)$] 
When $u$ converges to $u_\infty,$ we denote by $\vphi(s)=cs^\theta$ $\left(\text{with } c>0 \text{ and } \theta\in\left(0,\frac12\right]\right)$ the desingularizing function of $G$ at $u_\infty$ -- the quantity $\theta$ is the {\L}ojasiewicz exponent associated with $u_\infty.$ One has the following estimates.
\begin{itemize}
\item[$(a)$]
$\|u(t)-u_\infty\|\le ct^{-\frac{\theta}{1-2\theta}},$ with $c>0,$ when $\theta\in\left(0,\frac12\right).$
\item[$(b)$]
$\|u(t)-u_\infty\|\le c^\pp\exp(-c^\p t),$ with $c^\p,c^\pp>0,$ when $\theta=\frac12.$
\end{itemize}
\end{itemize}
\end{cor}

\begin{proof*}
The proof follows directly from the original {\L}ojasiewicz inequality \cite{MR0160856,Lo65} and the fact that desingularizing functions for real-analytic functions are indeed of the form $\vphi(s)=cs^\theta$ with $\theta\in(0,\frac12].$ Hence \eqref{grow_ineq} holds and Theorem \ref{thmmain} applies, see also Remark~\ref{vitesse} (c).
\medskip
\end{proof*}

\begin{cor}[\textbf{Convergence theorem for definable functions}]
Let $\vO$ be an o-minimal structure that contains the collection of semi-algebraic sets. Assume $G:\R^N\tends\R$ is $C^2$ and definable in $\vO.$ Let $u$ be a bounded solution to \eqref{second1}.
Then we have the following result.
\begin{itemize}
\item[$(i)$]
$u^\p$ and $u^{\p\p}$ belong to $L^1\big((0,\infty);\R^N\big)$ and in particular $u$ converges to a single limit $u_\infty$ in $\crit G.$
\item[$(ii)$]
When $u$ converges to $u_\infty$ we denote by $\vphi$ the desingularizing function of $G$ at $u_\infty.$ One has the following estimate
\begin{gather*}
\|u(t)-u_\infty\|\le c\nu(t),
\end{gather*}
where $\nu$ is a solution of the worst-case gradient system 
\begin{gather*}
\nu^\p(t)+(\vphi^{-1})^\p(\nu(t))=0, \; \nu(0)>0.
\end{gather*}
\end{itemize}
\end{cor}

\begin{proof*}
$G$ is a KL function by Kurdyka's version of the {\L}ojasiewicz inequality. The fact that $\vphi^\p(s)\ge\frac{c}{\sqrt{s}}$ comes from Lemma~\ref{lower}. So, Theorem~\ref{thmmain} applies.
\medskip
\end{proof*}

\begin{cor}[\textbf{Convergence theorem for the one-dimensional case \cite{MR869543}}]
Let $G\in C^{2}(\R;\R)$ and let $u$ be a bounded solution to~\eqref{second1}. Then $u$ converges to a single point and we have the same type of rate of convergence as in the previous corollary.
\end{cor}

\begin{proof*}
We proceed as in \cite{MR0223758}. Argue by contradiction and assume that $\o(u_0,u_0^\p),$ the $\o$-limit set of $(u_0,u_0^\p),$ is not a singleton. Since $\o(u_0,u_0^\p)$ is connected in $\R,$ it is an interval and has a nonempty interior. Take $\ovl u$ in the interior of $\o(u_0,u_0^\p)$ The {\L}ojasiewicz inequality trivially holds at $\ovl u$ for $G\equiv0$ with $\vphi(s)=\sqrt{s}$ (recall $\ovl u$ is interior). Apply then Theorem~\ref{thmmain}.
\medskip
\end{proof*}

\begin{rmk} In the one-dimensional case, convergence can be obtained with much more general forms of damping, see \cite{MR2529922}.
\end{rmk}

\begin{cor}[\textbf{Convergence theorem for convex functions satisfying growth conditions}]
\label{convex}
\text{} \\
Let $G\in C^2(\R^N;\R)$ be a convex function such that
\begin{gather*}
\argmin G\stackrel{\text{def}}{=}\big\{u\in \R^N;G(u)=\min G \big\},
\end{gather*}
is nonempty $($note that $\argmin G=\crit G).$ Assume further that, for each minimizer $x^*,$ there exists $\eta>0,$ such that $G$ satisfies
\begin{gather}
\label{grow}
\forall u\in B(x^*,\eta), \; G(u)\ge \min G + c\,\dist(u,\argmin G)^r, \, 
\end{gather}
with $r\ge1$ and $c>0.$ Then the solution curve $t\longmapsto(u(t),u^\p(t))$ has a finite length. In particular $u$ converges to a minimizer $u_\infty$ of $G$ as $t$ goes to $\infty.$
\end{cor}

\begin{proof*} A general result of Alvarez \cite{MR1760062} ensures that $u$ is bounded (and even converges).
On the other hand it has been shown in \cite{MR2592958} that functions satisfying the growth assumption \eqref{grow}, also satisfy the {\L}ojasiewicz inequality with desingularizing functions of the form $s\longmapsto c^\p s^{1-1/r}$ with $c^\p>0.$ Combining the previous arguments, the conclusion follows readily.
\medskip
\end{proof*}

\begin{rmk}
An alternative and more general approach to establish that trajectories have a finite length has been developed for convex functions in \cite{MR1104346,MR2665417}.
\end{rmk}

\appendix
\section{Appendix: some elements on o-minimal structures}
\label{appendix}

Some references for o-minimal structures are \cite{Coste99,MR1404337,MR1644089,MR1633348}. We only collect in this appendix the elements that are necessary to follow our main developments.

\begin{defi}[\textbf{o-minimal structure \cite[Definition\,1.5]{Coste99}}]
\label{domin}
An \textit{o-minimal} structure on $(\R,+,\:.\:)$ is a sequence of Boolean algebras$^($\footnote{Recall that a Boolean algebra is stable by finite union, finite intersection and contains the empty set and the total space; here $\emptyset\in\vO_n$ and $\R^n\in\vO_n.$}$^)$ $\vO=\{\vO_n\}_{n\in\N}$ of subsets of $\R^n$ such that for each $n\in\N,$
\begin{enumerate}
\itemsep=1mm
\item[(i)]
if $A$ belongs to $\vO_n$ then $A\times\R$ and $\R\times A$ belong to $\vO_{n+1};$
\item[(ii)]
if $\Pi:\R^{n+1}\tends\R^n$ is the canonical projection onto $\R^n$ then for any $A\in\vO _{n+1},$ the set $\Pi(A)$ belongs to $\vO_n;$
\item[(iii)]
$\vO_n$ contains the family of real algebraic subsets of $\R^n,$ that is, every set of the form
\begin{gather*}
\Big\{x\in\R^n; \, p(x)=0\Big\},
\end{gather*}
where $p:\R^n\tends\R$ is a real polynomial function;
\item[(iv)]
the elements of $\vO_1$ are exactly the finite unions of intervals and points.
\end{enumerate}
\end{defi}

\noindent
Being given an o-minimal structure $\vO,$ a set $A\subset\R^n$ is called \textit{definable} (in $\vO)$ if $A\in\vO_n.$ A mapping $F:D\subset\R^n\tends\R^m$ is said to be \textit{definable in} $\vO$ if its graph is definable in $\O$ as a subset of 
$\R^n\times\R^m.$ A point-to-set mapping
\begin{gather*}
S:\R^n\rightrightarrows\R^m,
\end{gather*}
maps each point $x$ in $\R^n$ to a subset $S(x)$ of $\R^m.$ The \textit{domain of} $S,$ denoted by $\dom S,$ is given by the set of elements $x$ in $\R^n$ such that $S(x)$ is nonempty. The graph of $S$ is defined by 
\begin{gather*}
\mathrm{graph}\,S=\Big\{(x,y)\in \R^n\times \R^m; y\in S(x) \Big\}.
\end{gather*}
As previously a point-to-set mapping is called \textit{definable $($in} $\vO)$ if its graph is definable in $\R^n\times\R^m.$ 

\begin{exa}
\label{exomin}
(a) \textbf{Semi-algebraic sets.}
The first and simplest example of o-minimal structure is given by the class of semi-algebraic objects (see \eqref{SA}). Tarski-Seidenberg principle (see \cite{MR1659509}) asserts that linear projections of semi-algebraic sets are semi-algebraic sets, in other words item (ii) of Definition~\ref{domin} holds for the class of semi-algebraic sets. The other items of the definition are easy to establish. \\
(b) \textbf{Globally subanalytic sets.}
There exists an o-minimal structure that contains semi-algebraic sets and sets of the form $\big\{(x,t)\in [-1,1]^n\times\R; \, f(x)=t\big\},$ where $f:[-1,1]^n\tends\R$ $(n\in \N)$ is a real analytic function that can be extended analytically on a neighborhood of the square~$[-1,1]^n$ -- these are sometimes called restricted analytic functions. This result is essentially due to Gabrielov \cite{MR1389958}; sets belonging to this structure are called {\em globally subanalytic sets} (see \cite{MR1289495} and the references therein). \\
(c) \textbf{Log-exp structure.}
There exists an o-minimal structure containing the globally subanalytic sets and the graph of $\exp:\R\tends\R,$ see \cite{MR1289495}.
\medskip
\end{exa}

\noindent
There are other results on o-minimal structures and the field is still very active, but the above examples give a good idea of the power of the concept.
\medskip \\
We now describe some stability/regularity results that we used in this paper. 
\medskip \\
Let $\vO$ be an o-minimal structure on $(\R,+,\:.\:).$

\begin{lem}[{\bf Monotonicity Lemma \cite[Theorem\,4.1]{MR1404337}}]
\label{mono}
Let $f:I\subset\R\tends\R$ be a definable function and $k\in \N.$ Then there exists a finite partition of $I$ into $p$ intervals $I_1,\ldots,I_p,$ such that $f$ restricted to each nontrivial interval $I_j,$ $j\in\{1,\ldots,p\},$ is $C^k$ and either strictly monotone or constant. Observe that some $I_j$ can be reduced to a singleton.
\end{lem}

\begin{lem}[{\bf Definable Selection Lemma \cite{Coste99}}]
\label{selec}
Let $S:\R^n\tends\R^m$ be a definable point-to-set mapping. Then there exists a definable mapping $F:\dom S\tends\R^m$ such that
\begin{gather*}
F(x)\in S(x),\: \forall x\in\dom S.
\end{gather*}
\end{lem}
\noindent
We recall the following theorem as stated in Kurdyka's original work \cite{MR1644089}.

\begin{thm}\label{KL0}
Let $\O$ be a nonempty open bounded subset of $\R^n$ and $f:\O \to\R$ a differentiable definable function with $f>0$ on $\O.$ Then there exist $r_0>0$ and a continuous definable function $\varphi:[0,r_0)\to \R_+$ such that $\varphi(0)=0,$ $\varphi\in C^1(0,r_0)$ and $\varphi^\p>0$ such that
\begin{gather*}
\|\nabla \left(\varphi\circ f\right)(x)\|\ge1,\: \forall x\in \O.
\end{gather*}
\end{thm}

\begin{rmk}
\label{KLK}
Let us show how to recover the form of KL inequality given in Theorem~\ref{KL}. 
\\
We adopt the notation of Theorem~\ref{KL}. Fix $\mu>0.$ Apply first, the above result to $G-G(\ovl u)$ (respectively, to $G(\ovl u)-G)$ on $\O_1=B(\ovl u,\mu)\cap[G-G(\ovl u)>0]$ (respectively, on $\O_2=B(\ovl u,\mu)\cap[G(\ovl u)-G>0]).$ This gives $\varphi_1:[0,r_1)\tends\R_+$ and $\varphi_2:[0,r_2)\tends\R_+,$ as in Kurdyka's Theorem. Let us now build a ``global" $\varphi$ as in Theorem~\ref{KL}. First recall that the derivative of a differentiable definable function is definable in the same structure, see \cite{Coste99}. Set $p(s)=(\varphi^\p_1-\varphi^\p_2)(s).$ By definability, $p$ is positive, negative or null on an interval of the form $(0,\eps).$ This yields the existence of $r$ in $(0,\min\{ r_1,r_2\})$ such that, for instance, $\varphi^\p_1>\varphi^\p_2$ on $(0,r).$ Set then $\varphi=\varphi_1$ and observe that
\begin{gather*}
\left\|\nabla\left(\varphi\circ|G(\:\cdot\:)-G(\ovl u)|\right)(u)\right\|\ge1,\: \forall u\in B(0,\eta)\setminus [G\neq G(\ovl u)],
\end{gather*}
when $\eta$ is sufficiently small.
\end{rmk}

\medskip

\noindent
\textbf{Acknowledgements.}
We are grateful to the referees for their very careful reading and their constructive input. 


\def\cprime{$^\p$}

\end{document}